\documentstyle[12pt]{article} %Remove [12pt] to print in 10pt

\title{On the space of harmonic $2$-spheres in ${\bf C}P^2$}
\author{L. LEMAIRE \thanks{C.P. 218 Campus
Plaine, Universit\'e libre de Bruxelles,
Bd.\ du Triomphe, 1050 Bruxelles, Belgium, e-mail: ulbmath@ulb.ac.be} $\,$
and J.C.
WOOD \thanks{Department of Pure Mathematics, University of Leeds, Leeds LS2
9JT, G.B., e-mail:
j.c.wood@leeds.ac.uk}} \date{ }

\newcommand{\Bbb}{\bf}   %Comment out to get \Bbb in place of \bf
\newcommand{{\CC}}{{\Bbb C}}
\newcommand{\RR}{{\Bbb R}}
\newcommand{\Hol}{\mbox{\rm Hol}}
\newcommand{\Harm}{\mbox{\rm Harm}}

\newtheorem{theo+}           {Theorem}      [section]
\newtheorem{prop+}  [theo+]  {Proposition}
\newtheorem{coro+}  [theo+]  {Corollary}
\newtheorem{lemm+}  [theo+]  {Lemma}
\newtheorem{exam+}  [theo+]  {Example}
\newtheorem{rema+}  [theo+]  {Remark}
\newtheorem{defi+}  [theo+]  {Definition}

\newtheorem{exam+s}  [theo+]  {Examples}
\newtheorem{rema+s}  [theo+]  {Remarks}
\newtheorem{hyp+}  [theo+]  {Hypotheses}
\newtheorem{cla+}  [theo+]  {Claim}

\newenvironment{theorem}{\begin{theo+}}{\end{theo+}}
\newenvironment{proposition}{\begin{prop+}}{\end{prop+}}

\newenvironment{lemma}{\begin{lemm+}}{\end{lemm+}}
\newenvironment{example}{\begin{exam+}\rm}{\end{exam+}}
\newenvironment{remark}{\begin{rema+}\rm}{\end{rema+}}
\newenvironment{definition}{\begin{defi+}\rm}{\end{defi+}}

\begin{document}
\maketitle
\begin{abstract}
Carrying further work of T.A. Crawford, we show that each
component of the space of harmonic maps from  the $2$-sphere to
complex projective $2$-space of degree $d$ and energy $4 \pi E$
is a smooth closed submanifold of the space of all $C^j$ maps $(j
\geq 2)$.  We achieve this by showing
that the Gauss transform which relates them to spaces of
holomorphic maps of given degree and ramification index is {\bf
smooth} and has {\bf injective differential}.
\end{abstract}

\section{Introduction}

In \cite{Cra}, T.A. Crawford showed that if $E \leq 5 |d|+10$, the
subset $\Harm_{d,E} (\CC P^2)$ of the space of harmonic maps from
$S^2$ to $\CC P^2$ consisting of those maps of degree $d$ and
energy $4\pi E$ can be given the structure of a complex manifold
and that this manifold is connected.

This he did by showing that the space $\Hol_{k,r}^*(\CC P^2)$ of
full holomorphic maps from $S^2$ to $\CC P^2$ of degree $k$ and
ramification index $r$ is
a complex manifold if $r \leq (k+1)/2$, and that the ``Gauss
transform'' $G'_{k,r}$ which maps $\Hol_{k,r}^*(\CC P^2)$ to
$\Harm_{k-r-2,3k-r-2}(\CC P^2)$ bijectively is a homeomorphism,
so that the manifold structure of $\Hol_{k,r}^*(\CC P^2)$ can be
transported to the topological space $\Harm_{d,E}(\CC P^2)$. This
does not prove that the transported structure is the one induced
by the natural inclusion of $\Harm_{d,E}(\CC P^2)$ in the space
of maps from $S^2$ to $\CC P^2$.

In this paper, after giving a treatment of Crawford's result
adapted to our needs, we show that $G'_{k,r}$ is a {\bf smooth}
map from $\Hol_{k,r}^*(\CC P^2)$ to $C^j(S^2,\CC P^2)$ \ (for $j
\geq 2$) and has {\bf injective differential}. From this we
obtain: \begin{theorem} \label{th:1.1}
For $0 \leq r \leq \displaystyle{\frac{k+1}{2}}$
and $\displaystyle{\frac{4k-11}{3} \leq r \leq
\frac{3}{2}k-3}$ \ the map
$$
G'_{k,r} : \Hol_{k,r}^* (\CC P^2) \to C^j(S^2,\CC P^2)
$$
is a smooth embedding onto $\Harm_{k-r-2,3k-r-2}(\CC P^2)$
for any $j \geq 2$.

Each component $\Harm_{d,E}(\CC P^2)$ of $\Harm(\CC P^2)$ with $E
\leq 5 |d| + 10$ is a closed  smooth submanifold of $C^j(S^2,\CC
P^2)$ of dimension $6E+4$ if $E=|d|$ (in which case it consists
of holomorphic or antiholomorphic maps) and of dimension $2E+8$
(otherwise). \end{theorem}

Added in proof: In a revised version of [4], Crawford has shown that
$\Hol_{k,r}^* (\CC P^2)$ is a manifold for
all $k,r$.  It follows that all our results (Theorems \ref{th:1.1},
\ref{th:1.3}, Proposition \ref{prop:3.1}, \ref{prop:4.1}, \ref{prop:5.2} and
Lemma \ref{lem:5.1}) are valid for this range and the restriction $E \leq 5 |d|
+ 10$ can be removed from Theorem \ref{th:1.1}.  Proofs are
unchanged except for those in Sec. 3, see below.

\begin{remark} \label{rem:1.2} The space $\Hol_{k,r}^*(\CC P^2)$
is non empty precisely for $k \geq 2$, \\ $0 \leq r \leq
\frac{3}{2} k-3$ (see Proposition \ref{prop:2.7} below).
\end{remark}

We shall first of all prove that $\Hol^*_{k,r}(\CC P^2)$ is a
complex manifold for the range $k \geq 2$, \ $0 \leq r \leq
(k+1)/2$\ . The passage from this range to the second $k \geq 3$,
\ $(4k-11)/3 \leq r \leq 3k/2 -3$ is achieved by the conjugate
polar (see Definition \ref{def:2.3} below).

In fact, we have:

\begin{theorem} \label{th:1.3} For $0 \leq r \leq
\displaystyle{\frac{k+1}{2}}$ and for
$\displaystyle{\frac{4k-11}{3} \leq r \leq \frac{3}{2}k-3}$,
$\Hol_{k,r}^*(\CC P^2)$ is a complex submanifold of the complex
manifold $\Hol_{k}^*(\CC P^2)$ of dimension $3k-r+2$.
\end{theorem}

Added in proof: Note that  Lemma \ref{lem:3.3} in the proof is
  false outside the above range; we give a counter example for $k=6$, $r=4$.
Hence the description of the manifold structure on $\Hol_{k,r}^*(\CC P^2)$
in Sec. 3 is special
to the range $r \leq (k+1)/2$ --- see the revised version
of \cite{Cra} for a description valid for any $k,r$.

\bigskip

The contents of the subsequent sections are as follows.

In Section 2, we recall the construction of J. Eells and the
second author of the harmonic maps of $S^2$ to $\CC P^2$,
stressing the limits of the values of the parameters involved,
and we give some examples illustrating the behaviour of this
construction.

In Section 3, we present a proof of the result of T.A. Crawford
on the structure of $\Hol_{k,r}^*(\CC P^2)$, adapted to our
needs. We use ideas from the paper of Crawford, and a
construction suggested to us by M. Guest.

In Sections 4 and 5, we show successively that the Gauss
transform is smooth and that it is an embedding. The main
ingredient is a property of smooth dependence of the common
factor of smooth families of polynomials (Lemma \ref{lem:4.5}).
We note at this point that the methods of the present paper will
not generalize easily to ${\CC} P^n$ for $n > 2$, however M.
Guest and Y. Ohnita \cite{G-O} show that some topological
questions on $\CC P^n$ reduce to $\CC P^2$.

Recall finally that the harmonic maps from $S^2$ to any manifold
are precisely the minimal branched immersions of $S^2$, in
particular, they are conformal (see e.g. \cite{E-L} (5.15))

\noindent
{\bf Acknowledgements}\\
Both authors would like to thank the organizers of the first MSJ
International Research Institute (Sendai, 1993) for invitations
to the conference where this work was begun, F. Burstall, J.
Eells, M. Guest and A. West for helpful discussions, and L. Simon
for posing questions at the above conference which led to this
research.

The second author would like to thank the Belgian Contact Group
for an invitation to their 1993 meeting which allowed this work
to continue.

This research was partially supported by the EU Human Capital and
Mobility Programme (Contract CHRX-CT92-0050), the Belgian
F.N.R.S. and the British Council.

\section{Gauss transform and polars}

Following work of others \cite{G-S,D-Z}, see also \cite{Bur}, J.
Eells and the second author classified the harmonic maps from the
2-sphere $S^2$ to a complex projective space ${\CC} P^n$, as
follows:

\begin{theorem} \label{th:2.1} \cite{E-W}. There is a bijective
correspondence between the set of pairs $(f,s)$ where $f : S^2
\to C P^n$ is a full holomorphic map and $s$ an integer, $0 \leq
s \leq n$, and the set of full harmonic maps $\varphi : S^2 \to
\CC P^n$.
\end{theorem}

Here ``full'' means ``not having image in any proper projective
subspace of $\CC P^n$\ ''. For $s=0$, $\varphi=f$ is holomorphic,
for $s \neq 0,n$, $\varphi$ is neither holomorphic nor
antiholomorphic, and for $s=n$, $\varphi$ is antiholomorphic and
is called the {\bf polar} of $f$.

We now restrict to the case of $\CC P^2$, where the description
of this construction is simpler.

Let $f : S^2 \to \CC P^2$ be a holomorphic map. Identifying $S^2$
with $\CC \cup \{ \infty \}$ by stereographic projection, $f$ can
be represented on $\CC$ by a map $p : \CC \to \CC^3 \setminus
\{0\}$ where $p(z) = (p_0(z),p_1(z),p_2(z))$ is a triple of
coprime polynomials with max (degree $p_0$, degree $p_1$, degree
$p_2$) = degree of $f$. We shall write $f=[p_0,p_1,p_2]$.

A map $S^2 \to \CC P^2$ is called {\bf full} if its image lies in
no complex projective line. Note that if a harmonic map is not
full, its image lies in a $\CC P^1$and it is then
$\pm$-holomorphic, since it is a conformal map between surfaces.
Thus all harmonic non $\pm$-holomorphic maps $f : S^2 \to \CC
P^2$ are full.

We denote by $\Hol_k^*(\CC P^2)$ the space of full holomorphic
maps of degree $k$. All values of $k \geq 2$ occur and $\Hol_k^*
(\CC P^2)$ is a complex manifold of dimension $3k+2$ with
coordinate charts given by the coefficients.

Recall that a holomorphic map $f : S^2 \to \CC P^2$ is said to be
{\bf ramified} at a point $z \in S^2$ if $df(z)=0$. The
ramification index of $f$ at $z$ is the order of the zero of
$df(z)$ and the {\bf total ramification index} $r$ of $f$ is the
sum of ramification indices.

Consider a full holomorphic map $f$. The harmonic map $\varphi$
associated to the pair $(f,1)$ in Theorem \ref{th:2.1} is
obtained by means of the $\partial'$-{\bf Gauss transform} (in
the terminology of \cite{B-W} - in \cite{Wol} it is called the
$\partial$-{\bf transform}) which is defined as follows
\cite{E-W}:

Let $\pi : \CC^3 \setminus \{0\} \to \CC P^2$ be the canonical
projection sending $(z_0,z_1,z_2)$ to the point of $\CC P^2$ with
homogeneous coordinates $[z_0,z_1,z_2]$. For a map $f : S^2 \to
\CC P^2$, say that $F : U \to \CC^3 \setminus \{0\}$ represents
$f$ on the open set $U$ if $f |_U = \pi \circ F$, in which case
we write $f=[F]$.

For $f$ holomorphic and full, the {\bf first associated curve}
$f_{(1)} : S^2 \to G_2(\CC^3)$ is the holomorphic map defined as
follows. Let $F : U \to \CC^3 \setminus \{0\}$  represent $f$ on
a domain $(U,z)$ of $\CC$. Consider the map $F \wedge F' : U \to
\wedge^2 \CC^3$ where $\displaystyle{F'=\frac{dF}{dz}}$. At a
point $x$ where $f$ is not ramified, $F \wedge F'$ is non zero
and so defines a complex two-dimensional subspace $f_{(1)}(x)$.
If, on the other hand, $f$ is ramified at $x$ with ramification
index $k$, then $F \wedge F' = (z-x)^k \cdot \Psi$ for some
smooth nonzero map $\Psi : U' \to \wedge^2 \CC^3$ on an open
neighbourhood of $x$. Since $\Psi(y)$ is decomposable for all $y
\in U'$, $y \neq x$, it remains decomposable for $y=x$ and we can
define $f_{(1)}(x)$ as the complex two-dimensional subspace
defined by $\Psi(x)$. The resulting map $f_{(1)} : S^2 \to
G_2(\CC^3)$ is well-defined and smooth. This leads us to two maps
associated to $f$~:

\begin{definition} \label{def:2.2} The {\bf $\partial'$-Gauss
transform} $\varphi=G'(f) : S^2 \to \CC P^2$ is defined by the
formula \begin{equation}\label{equ:2.1} \varphi(x) = f(x)^\perp
\cap f_{(1)}(x). \end{equation} \end{definition}

\noindent By Theorem 2.1, it is a smooth and full harmonic map.

\begin{definition} \label{def:2.3} The {\bf polar} of the holomorphic
map $f$ is the antiholomorphic map $$ g(x) =f_{(1)}(x)^\perp \ .
$$
\end{definition}

Note that for any $x \in S^2$, $f(x), \varphi(x)$ and $g(x)$ are
Hermitian orthogonal complex lines.

For convenience, we consider the {\bf conjugate polar} $h$ of $f$
defined by taking in $\CC^3$ the complex conjugate of the values
of $g(x) : h(x) = \overline{f_{(1)}(x)^\perp}$.

More explicitly, represent $f \in \Hol_k^*(\CC P^2)$ by
$[p_0,p_1,p_2]$ as above.  Identifying $\wedge^2 \CC^3$ with
$\CC^3$, the first associated curve, or equivalently the
conjugate polar $h$ of $f$, is represented by the polynomials
\begin{eqnarray*} h &=& [p_{12},p_{20},p_{01}]\\ &=&
[p_1p'_2-p'_1p_2,p_2p'_0-p'_2p_0,p_0p'_1-p'_0p_1] \end{eqnarray*}
once they have been divided by their common factor. Explicitly,
if $f$ has finite ramification points $z_I$ with multiplicities
$r_I$ $(I=1,\ldots,R)$, then the {\bf ramification divisor}
$R(f)$ is the monic polynomial $R(f)=\prod (z-z_I)^{r_I}$.

If $f$ is not ramified at $\infty$ this has degree equal to the
ramification index, otherwise it has lower degree.

In either case, $R(f)$ is the highest common factor of the
polynomials $p_{ij}$, and $$ h = \left[ \frac{p_{12}}{R(f)},
\frac{p_{20}}{R(f)}, \frac{p_{01}}{R(f)} \right]. $$ One checks
easily that $h$ has degree $2k-2-r$ (indeed, the terms of degree
$2k-1$ in $p_{ij}$ cancel).

Note that formula (\ref{equ:2.1}) has a counterpart:
$$
\overline{\varphi(x)} = h(x)^\perp \cap h_{(1)} (x).
$$
Theorem \ref{th:2.1} can be rephrased in $\CC P^2$ by saying that
{\em the Gauss transform defines a bijection between the space of
full holomorphic maps and the space of harmonic non
$\pm$-holomorphic maps}, and that {\em the passage to the polar
defines a bijection between full holomorphic and antiholomorphic
maps}.

However, this does not immediately provide a simple description
of the space of harmonic maps. Indeed, the Gauss transform $G':
\Hol_k^*(\CC P^2) \to \Harm(\CC P^2)$ is not a continuous map,
when the spaces are equipped with their $C^0$-topology. This
appears for instance in the following example (brought to our
attention by F. Burstall).

\begin{example} \label{ex:2.4} Let $f_t : S^2 \to \CC P^2$ be
defined by $f_t(z) = [F_t(z)]$, where $$ F_t(z) = (1,tz+z^3,z^2)
\quad (z \in \CC, t\in {\bf R}) $$ (so that
$f_t(\infty)=[0,1,0]$). Note that $f_t(0)=[1,0,0]$ for all $t$
and that $f_t$ is a smooth family of full holomorphic maps. Then
$F'_t(z)=(0,t+3z^2,2z)$ and $F_t \wedge
F'_t(z)=(tz^2-z^4,-2z,t+3z^2)$. If $t \neq 0$, at $z=0$ this
equals $(0,0,t)=t(1,0,0) \wedge (0,1,0)$ so that the first
associated curve has the value $f_{t(1)}(0) =
\mbox{span}\{(1,0,0),(0,1,0)\}$ and $G'(f_t)(0)=[0,1,0]$.
However, if $t=0$, then $F_t \wedge
F'_t(z)=(-z^4,-2z,3z^2)=z\psi(z)$ where $\psi(z)=(-z^3,-2,3z)$.
In particular $\psi(0)=(0,-2,0)=2(1,0,0) \wedge(0,0,1)$ so that
$f_{0(1)}(0) = \mbox{span}\{(1,0,0),(0,0,1)\}$ and
$G'(f_0)(0)=[0,0,1]$. This shows that $f_{t(1)}$ and $G'(f_t)$ do
not vary continuously with $t$. The reason for this is that $f_t$
is unramified when $t \neq 0$ but ramified with ramification
index 1 at $z=0$ when $t=0$, and that in the presence of
ramification, both $G'$ and $g$ involve division of polynomials
by their common factor, a discontinuous process when the degree
of the factor changes.
\end{example}

To proceed, define $\Hol_{k,r}^*(\CC P^2)$ as the space of full
holomorphic maps of degree $k$ and total ramification index $r$,
and $\Harm_{d,E}(\CC P^2)$ the space of all harmonic maps of
degree $d$ and energy $4 \pi E$. By results of \cite{E-W},
Theorem \ref{th:2.1} specializes to

\begin{proposition} \label{prop:2.5} For each pair of integers $k
\geq 2$, $0 \leq r \leq \frac{3}{2}k-3$, there is a bijective
correspondence $$ G'_{k,r} : \Hol_{k,r}^*(\CC P^2) \to
\Harm_{d,E}(\CC P^2) $$ given by the restriction of the Gauss
transform, where $d = k-r-2$ and $E=3k-r-2$.
\end{proposition}

\begin{proposition} \label{prop:2.6} For each pair of integers $k
\geq 2$, $0\leq r \leq \frac{3}{2}k-3$, the map $f \mapsto$
conjugate polar of $f$ restricts to a bijection $$
\Hol_{k,r}^*(\CC P^2) \to \Hol_{k',r'}^*(\CC P^2) $$ where
$k'=2k-r-2$, $r'=3k-2r-6$.
\end{proposition}

This allows us to specify the values of $k$ and $r$ as follows:

\begin{proposition} \label{prop:2.7} The space $\Hol_{k,r}^*(\CC
P^2)$ is non-empty precisely for the range $k \geq 2$, $0 \leq r
\leq \frac{3}{2}k-3$. \end{proposition}

\noindent
{\bf Proof.} The well known examples (see \cite{C-M-R}) $$ f(z) =
[1,(z+1)^{k-r+1},z^k] $$ provide maps $f$ in $\Hol_{k,r}^*(\CC
P^2)$ for all $k \geq 2$, $0 \leq r \leq k-2$.

The Pl\"ucker formulae (see e.g. \cite{E-W}) show that the
involutive map $f \mapsto$ conjugate polar of $f$ restricts to
bijections $$ \Hol_{k,r}^*(\CC P^2) \to \Hol_{k',r'}^*(\CC P^2)
$$ with $k',r'$ as in Proposition 2.6. Thus to get
$\Hol_{k,r}^*(\CC P^2)$ (and $\Hol_{k',r'}^*(\CC P^2)$) non
empty, we need $r' \geq 0$, i.e. $r \leq \frac{3}{2}k-3$.

On the other hand, the conjugate polars of the above examples
provide maps in $\Hol_{k,r}^*(\CC P^2)$ for all $$ k-2 \leq r
\leq \frac{3}{2} k-3 \ . $$ \medskip

As in \cite{C-M-R}, we note that the above shows that
$\Harm_{d,E} (\CC P^2)$ is non-empty precisely for pairs $(d,E)$
of integers with either $E=|d|$ (in which case it consists of
$\pm$-holomorphic maps) or $E=3|d|+4+2r$ for some $r \geq 0$
(otherwise).

Indeed, all such values of $E$ with $d \geq 0$ are achieved with
$0 \leq r \leq k-2$, the range $k-2 < r \leq 3k/2-3$ giving $d <
0$.

The main contribution of the present paper is to show that $$
G'_{k,r} : \Hol_{k,r}^*(\CC P^2) \to \Harm_{d,E}(\CC P^2) \subset
C^j(S^2,\CC P^2) $$ is a smooth embedding and that the map $f
\mapsto$ conjugate polar of $f$ is a complex analytic equivalence
from $\Hol_{k,r}^*(\CC P^2)$ to $\Hol_{k',r'}^*(\CC P^2)$.

We conclude this section by an example showing that in
$\Hol_{k,r}^*(\CC P^2)$, the ramification divisor of a smooth
family of polynomials can vary smoothly, even when individual
common roots of the $p_{ij}$'s vary only continuously.

\begin{example} \label{ex:2.8} Identifying $S^2$ with $\CC \cup
\{0\}$ by stereographic projection, let $f_t : S^2 \to \CC P^2$
be defined by $$ F_t(z)=\big
(z^4+1,(1-3t^2)z^3+(-3t+t^3)z,2tz^2+(1-t^2) \big) \quad (z \in
\CC, t \in {\RR}) $$ (so that $f_t(\infty)=[1,0,0]$).

Identifying $\Lambda^2 \CC^3$ with $\CC^3$ we have $(F_t \wedge
F'_t)(z)=(z^2-t)\psi(z)$ where $$ \psi(z) = \big(
(-2t+6t^3)z^2+(-3+t^2)(1-t^2),4z(tz^2+1),(-1+3t^2)z^4+8tz^2+3-t^2
\big) \ . $$ which shows that, if $t \neq 0$, $f_t$ is ramified at
$z=\pm \sqrt{t}$ with index 1, but if $t=0$, these ramification
points coalesce into a ramification point at $z=0$ of index 2.
Further $f_{t(1)}(z)=[\psi(z)]$. We see from this that $f_t \in
\Hol_{4,2}^*(\CC P^2)$ for all $t$ and that $f_{t(1)}$ and so
$G'(f_t)$ vary smoothly with $t$, even though each root does not.
\end{example}

\section{Spaces of holomorphic maps} In this section, we give a
proof of the

\begin{proposition} \label{prop:3.1} \cite{Cra} For any $k \geq
2$ and $0 \leq r \leq \displaystyle{\frac{k+1}{2}}$,
$\Hol_{k,r}^* (\CC P^2)$ is a complex submanifold of
$\Hol_k^*(\CC P^2)$ of dimension $3k-r+2$.
\end{proposition}

\noindent {\bf Proof.} Let
\begin{eqnarray*}
\Hol'_k(\CC P^2) & = & \{ f \in \Hol_k^*(\CC P^2) : f =
[p_0,p_1,p_2],\ p_0 \mbox{ is monic of}\\ & & \mbox{degree } k
\mbox{ with distinct roots}, f \mbox{ is not ramified at } \infty
\} \ .
\end{eqnarray*}
Here, $p_0,p_1$ and $p_2$ are always assumed to be coprime.
$\Hol'_k(\CC P^2)$ is an open subset of $\Hol_k^*(\CC P^2)$ and
so a complex manifold. It can be embedded as an open subset in
$\CC^{3k+2}$ by sending the polynomials $(p_0,p_1,p_2)$ to their
coefficients (omitting the leading coefficient of $p_0$, which is
equal to 1).

The group $G=PGL_2(\CC) \times PGL_3(\CC)$ acts on $\Hol_k^*(\CC
P^2)$ in a natural way preserving the subsets $\Hol_{k,r}^*(\CC
P^2)$. Given $f \in \Hol_k^*(\CC P^2)$, a variation of the proof
of Lemma \ref{lem:4.5} below shows that there exist $g \in G$ and
a neighbourhood $U$ of $f$ in $\Hol_k^*({\bf C } P^2)$ such that
$g(U) \subset \Hol'_k(\CC P^2)$. Setting $\Hol'_{k,r}(\CC
P^2)=\Hol_{k,r}^*(\CC P^2) \cap \Hol'_k(\CC P^2)$, it suffices
therefore to show that $\Hol'_{k,r}(\CC P^2)$ is a complex
submanifold of $\Hol'_k(\CC P^2) \subset {\CC}^{3k+2}$.

To do this, we use a construction suggested by M. Guest following
ideas of T.A. Crawford. Set \begin{eqnarray*}
   X''_{k,r}  & = & \{ (a,f)=(a, (p_0,p_1,p_2)) \in \CC^r \times
\CC^{3k+2}~:\ f \in \Hol'_k(\CC P^2) \ , \\
               & &  a \mbox{ is a monic polynomial of degree } r ,
        (a,p_0)\mbox{ are coprime and } a  \mbox{ divides } R(f) \} \ .
\end{eqnarray*}
Note that $X''_{k,r}$ is an algebraic subvariety of $\CC^r \times {\bf
C}^{3k+2}$.
We shall prove that $X''_{k,r}$ is a complex submanifold provided
$r \leq (k+1)/2$.

There is an injective map
$$
i : \Hol'_{k,r}(\CC P^2) \to X''_{k,r}
$$
given by $i(f)=(R(f),(p_0,p_1,p_2))$ where $f=[p_0,p_1,p_2]$ as
above, with image $X'_{k,r}=\{(a,f) \in X''_{k,r} : \deg f = k,
\mbox{ ramification index } (f)=r \}$.

To check this, we need only show that in
$(R(f),f)=(a,(p_0,p_1,p_2))$, $a$ and $p_0$ are coprime.

Suppose, to the contrary, that there is an $x$ with
$a(x)=p_0(x)=0$. Then $p_0(z)=(z-x)p(z)$, and since $a$ divides
$p_0p'_i-p_ip'_0$, it follows that $p_i(x)p(x)=0$ for $i=1,2$.
Since $p_0$ has distinct roots, $p(x) \neq 0$ so that $p_i(x)=0$
for $i=1,2$, contradicting the fact that $p_0,p_1$ and $p_2$ are
coprime.

By Lemma \ref{lem:4.5} below, the injective map $i$ is complex
analytic, since the map $f \mapsto R(f)$ is.

Note that since $f \in \Hol'_{k,r}(\CC P^2)$ is not ramified at
infinity, $R(f)$ is a polynomial of degree $r$. The complement of
$X'_{k,r}$ in $X''_{k,r}$ is a proper subvariety of $X''_{k,r}$,
so that if $X''_{k,r}$ is a complex submanifold of $\CC^r \times
\CC^{3k+2}$, so is $X'_{k,r}$.

To study $X''_{k,r}$, we embed it in the trivial holomorphic
vector bundle $\pi : E \to A$, where $A$ is the open set in
$\CC^{r+k}$ given by \begin{eqnarray*} A & = & \{(a,p_0) \in
\CC^r \times \CC^k : a \mbox{ and } p_0 \mbox{ are monic coprime
polynomials } \\ & & \mbox{of degrees } r \mbox{ and } k \mbox{
respectively and } p_0 \mbox{ has no repeated root} \} \ , \\ E &
= & \{(a,p_0,p_1,p_2):(a,p_0)\in A, (p_1,p_2) \in \CC^{k+1}
\times \CC^{k+1} \} \end{eqnarray*} and $\pi$ is the natural
projection.

For $(a,p_0) \in A$, let $T_{(a,p_0)} : \CC^{k+1} \to \CC^r$ be
the linear map which sends a polynomial $p$ of degree $\leq k$
(represented by its coefficients) to the remainder of the
division of $p_0p'-p'_0p$ by $a$. We have:

\begin{lemma} \label{lem:3.2} Let $(a,p_0,p_1,p_2) \in E$ and
$f=[p_0,p_1,p_2]$. Then $a | R(f)$ if and only if $p_1$ and $p_2$
lie in $\ker T_{(a,p_0)}$.
\end{lemma}

\noindent
{\bf Proof.} With the notation $p_{ij}=p_ip'_j-p_jp'_i$, we have
immediately $$ p_1  p_{02} - p_2 p_{01}=p_0p_{12}. $$

If $p_1,p_2 \in \ker T_{(a,p_0)}$, then $a$ divides the left hand
side, and since $a$ and $p_0$ are coprime, $a$ must divide
$p_{12}$. Therefore $a$ divides $R(f)$.

The converse is immediate.

\medskip

Now, note that $X''_{k,r}$ is the kernel of the morphism of
holomorphic vector bundles
$$
E=A \times (\CC^{k+1})^2 \to A \times (\CC^r)^2
$$
defined by
$$
((a,p_0),(p_1,p_2)) \to ((a,p_0),T_{(a,p_0)}(p_1),T_{(a,p_0)}(p_2)).
$$
Hence $X''_{k,r}$ is a complex submanifold of $\CC^r \times
{\CC}^{k+2}$ if $\dim \ker T_{(a,p_0)}$ is independent of
$(a,p_0) \in A$.

\begin{lemma} \label{lem:3.3} If $\displaystyle{r \leq
\frac{k+1}{2}}$, $\dim \ker T_{(a,p_0)}=k+1-r$ \ $\forall (a,p_0)
\in A$. \end{lemma}

\noindent
{\bf Proof.} Let the zeros of $a$ be $\alpha_1, \ldots, \alpha_R$ with
multiplicities $m_1,\ldots,m_R$, so that $\displaystyle{\sum_{J=1}^R
m_J=r}$. For any $p$, set $h(p)=p_0p'-p'_0p$. Then $p \in \ker
T_{(a,p_0)}$ iff
\begin{equation}\label{equ:3.1}
(h(p))^{(I)}(\alpha_J)=0 \quad \forall J = 1,\ldots,R,\;
I=0,\ldots,m_J-1
\end{equation}
where $(h(p))^{(I)}$ denotes the $I^{\mbox{th}}$ derivative of the
polynomial $h(p)$.

Now (\ref{equ:3.1}) is a system of $r$ linear equations in $k+1$
unknows. Indeed, we can replace $T_{(a,p_0)}$ by the linear map
$\CC^{k+1} \to \CC^r$ which sends $p \in \CC^{k+1}$ to the vector
$$
\left( (h(p))^{(I)}(\alpha_J), \, J=1,\ldots,R, \, I=0,\ldots,m_J-1 \right)
\in \CC^r.
$$

We shall show that this map has rank $r$ by finding $r$
polynomials $P_{K,L} \in \CC^{k+1}$ $(L=1,\ldots,R$, $K=1,
\ldots,m_L)$ whose images are linearly independent.

To do this, choose for $P_{K,L}$ a polynomial of degree $\leq k$
with roots $\alpha_L$ of multiplicity $K$ and $\alpha_J$ (for $J
\neq L$) of multiplicity $m_J+1$.

This is possible since \begin{eqnarray*}
m_L + \sum_{J \neq L}(m_J+1) &=& r + R -1 \\
&\leq& 2r-1 \leq k
\end{eqnarray*}
by the hypothesis $r \leq (k+1)/2$.

\noindent
Then
$$
(h(P_{K,L}))^{(I)}(\alpha_J) = \left\{
\begin{array}{lll}
0 & &\mbox{if } J \neq L \\
&\mbox{or} &J=L \mbox{ and } I < K-1 \ . \\
\mbox{non zero } & &\mbox{if } J=L,\ I=K-1
\end{array}
\right.
$$
If we order the components of the vector $(h(p))^{(I)}(\alpha_J)$
in lexicographical order, viz.
$$
(J,I)=(1,0),(1,1),\ldots,(1,m_1-1),(2,0),\ldots,(R,m_R-1),
$$
we observe that the matrix
$$
(h(P_{K,L}))^{(I)}(\alpha_J)
$$
is in echelon form.

\noindent
Thus, the images of the $P_{K,L}$ are linearly independent, which
shows that the rank of $T_{(a,p_0)}$ is $r$, and so its kernel
has dimension $k+1-r$.

\begin{remark} \label{rem:3.4} An example, obtained in conjunction
with M. Guest, shows this lemma to be false for $k=6,r=4$.
Namely, when
$$
p_0=4z^6-12z^5+10z^4+2z^2-4z+4
$$
and
$$
a = z(z-1)(z+1)(z-2),
$$
the kernel of $T_{(a,p_0)}$ is of dimension 4, instead of 3.
\end{remark}

We deduce from the lemma that  $X''_{k,r}$, and so
$X'_{k,r}$, is a complex submanifold of dimension $3k-r+2$.

Now the restriction of the projection $X'_{k,r} \to \Hol'_k(\CC
P^2)$ which forgets $a$ is complex analytic and has image
$\Hol'_{k,r}(\CC P^2)$ (indeed, it is the inverse of the map $i :
\Hol'_{k,r}(\CC P^2) \to X''_{k,r}, \, f \mapsto (R(f),f))$.

Thus $\Hol'_{k,r}(\CC P^2)$ is a complex analytic submanifold of
$\CC^{3k+r+2}$, which concludes the proof of Proposition
\ref{prop:3.1}.

\section{The smooth nature of $G'_{k,r}$} In this section we
shall prove, with notation as in \S 1,

\begin{proposition} \label{prop:4.1} For $0 \leq r \leq
\displaystyle{\frac{k+1}{2}}$ and for
$\displaystyle{\frac{4k-11}{3} \leq r \leq \frac{3}{2}k-3}$, \
$G'_{k,r} : \Hol_{k,r}^*(\CC P^2) \to C^j (S^2,\CC P^2)$ is a
$C^\infty$ map between $C^\infty$ manifolds, for any $2 \leq j <
\infty$. \end{proposition}

\begin{lemma} \label{lem:4.2} Let $g_t$ and $h_t$ be two families
of polynomials in a single (complex) variable which depend
smoothly (resp. complex analytically) on a parameter $t \in U
\subset \RR^N$ (resp. $\CC^N$), where $U$ is open. Suppose that
the degrees of $g_t$, $h_t$ and of their highest common factor
$l_t$ are all constant, i.e. do not vary with $t$. Then the
polynomial $l_t$ depends smoothly (resp. complex analytically) on
$t$.
\end{lemma}

\begin{remark} \label{rem:4.3} We can take the polynomial $l_t$ to
be monic; then the statement of the lemma means that the
remaining coefficients depend smoothly (or complex analytically)
on $t$. \end{remark}

\begin{remark} \label{rem:4.4} In the situation of Lemma
\ref{lem:4.2}, each root of the polynomials depends continuously
on $t$, but not smoothly, in general, when the multiplicity
changes. So the linear factors of the polynomials do not always
vary smoothly. \end{remark}

\noindent
{\bf Proof of Lemma \ref{lem:4.2}.} It is sufficient to prove the
lemma for $t$ close to a fixed $t_0 \in U$. We can suppose deg
$l_t > 0$, or there is nothing to prove.

\noindent
{\bf First case:} Suppose that for one value of $t$, $g_t$ divides
$h_t$ or $h_t$ divides $g_t$. If for instance $g_t$ divides
$h_t$, we have $l_t=a(t)g_t$, where $a(t)$ is a scalar function.
Since all degrees are constant, we see that, for all $s$, deg
$l_s=$ deg $g_s$, so that again $l_s=a(s)g_s$, and $l_s$ is
smooth.\\
{\bf Second case:} For each $t,\ g_t$ (resp. $h_t$) does not divide
$h_t$ (resp. $g_t$). In particular, $g_t$ and $h_t$ are not
proportional.

For brevity of notation we shall now omit the parameter $t$ from
the notation.

{\bf Claim 1.} There exist unique polynomials $\lambda$ and $\mu$
with \begin{equation}\label{equ:4.1}
\mbox{deg } \lambda < \mbox{ deg } (h/l) \quad \mbox{ and } \quad
\mbox{ deg } \mu < \mbox{ deg } (g/l)
\end{equation}
such that $\lambda g + \mu h=l$.

The Euclidean algorithm ensures the existence of $\lambda$ and
$\mu$ such that $\lambda g + \mu h = l$. Suppose deg $\lambda
\geq$ deg $(h/l)$ and let $\tilde{\lambda}$ be the unique
polynomial such that
$$
\lambda = q \cdot \frac{h}{l} + \tilde{\lambda},
$$
with $\tilde{\lambda}=0$ or deg $\tilde{\lambda} <$ deg $(h/l)$.

If $\tilde{\lambda} = 0$, we have $\lambda = q \cdot h/l$ and
$\lambda g + \mu h = 1$ becomes $(qg/l+\mu)h=l$, which is
impossible with $l \not \equiv 0$ and deg $l <$ deg $h$.

If, instead, $\tilde{\lambda} \neq 0$, and deg $\tilde{\lambda}
<$ deg $(h/l)$, we have $\tilde{\lambda}=\lambda -q \cdot h/l$.
\\ Setting $\tilde{\mu} = \mu + q \cdot g/l$, we see that $$
\tilde{\lambda} g + \tilde{\mu}h = l
$$
and
$$
\mbox{deg } \tilde{\mu} + \mbox{ deg } h = \mbox{ deg }
\tilde{\lambda} + \mbox{ deg } g.
$$
This implies
\begin{eqnarray*}
\mbox{deg } \tilde{\mu} &<& \mbox{ deg } h - \mbox{ deg } l +
\mbox{ deg } g - \mbox{ deg } h \\ &=& \mbox{ deg } (g/l) \ .
\end{eqnarray*}

Thus we have existence of $\tilde{\lambda}$ and $\tilde{\mu}$
satisfying (\ref{equ:4.1}).

Unicity is easily checked.

{\bf Claim 2}. Let $\lambda$ and $\mu$ be polynomials satisfying
(\ref{equ:4.1}) and such that deg$(\lambda g + \mu h)  \leq$ deg
$l$. Then $\lambda g + \mu h=a \cdot l$, with $a(t)$ a scalar
function.

Since $l$ divides $g$ and $h$, we have deg$(\lambda g + \mu h)
\geq$ deg $l$ or $\lambda g + \mu h=0$.

In the second case, we have $\lambda g/l=- \mu h/l$, and by
(\ref{equ:4.1}) $g/l$ must have a common factor with $h/l$, a
contradiction.

Therefore, deg $(\lambda g + \mu h) \geq$ deg $l$. With the
hypotheses of the claim, we have deg $(\lambda g + \mu h)=$ deg
$l$ and $l$ divides $\lambda g + \mu h$ so that $\lambda g + \mu
h = a \cdot l$. \medskip

We conclude from the two claims that for $\lambda$ and $\mu$
satisfying (\ref{equ:4.1}), $l$ is characterized up to a non-zero
scalar factor as the unique polynomial of the form $\lambda g +
\mu h$ such that deg $(\lambda g + \mu h) \leq$ deg $l \equiv
L$.

Writing the parameter $t$ back in, this is equivalent to
$\displaystyle{\frac{d^{L+1}}{dz^{L+1}}(\lambda_tg_t+\mu_th_t)=0}$,
a system of homogeneous linear equations in the unknown
coefficients of $\lambda_t$ and $\mu_t$, with coefficients smooth
in $t$.

At $t=t_0$, consider any non zero coefficient of $\lambda_t$ or
$\mu_t$ and scale the solution by setting the coefficient equal
to 1 for $t$ close to $t_0$. The system becomes inhomogeneous,
and can be solved by Cramer's rule, so that the solutions are
smooth, which proves Lemma \ref{lem:4.2}.

\begin{lemma} \label{lem:4.5} Let $g_t, h_t$ and $k_t$ be three
families of polynomials in a single complex variable which depend
smoothly (resp. complex analytically) on a parameter $t \in U
\subset \RR^N$ (resp. $\CC^N$). Suppose that the degrees of $g_t$
and of the highest common factor $l_t$ of $g_t, h_t, k_t$ are
constant and that deg $h_t \leq$ deg $ g_t$ and deg $k_t \leq$
deg $g_t$ for all $t$. Then $l_t$ depends smoothly (resp. complex
analytically) on $t$.
\end{lemma}

\noindent
{\bf Proof.} The idea of the proof is to replace $g_t, h_t$ and
$k_t$ by linear combinations $\tilde{g}_t, \tilde{h}_t$ and
$\tilde{k}_t$, so that the common factor remains the same, but
any two of the three polynomials have no further common factor.

First, replace $h_t$ by $h_t + a \cdot g_t$ and $k_t$ by $k_t + b
\cdot g_t$ so that the three polynomials (still denoted by $g_t,
h_t$ and $k_t$) all have the same constant degree.

Consider now a fixed value of the parameter -- say $t=0$. We
shall show that there exists $\epsilon > 0$ such that $l_t$ is
smooth for $\| t \| < \epsilon$.

Let $A$ be a 3 by 3 matrix, which we shall choose close to the
identity matrix $I$, and in particular invertible. Set \begin{equation}
\left( \begin{array}{c}
\tilde{g}_t \\
\tilde{h}_t \\
\tilde{k}_t
\end{array}
\right)
 = A
\left( \begin{array}{c}
{g}_t \\
{h}_t \\
{k}_t
\end{array}
\right) \ . \label{equ:4.2}
\end{equation}

For $\| t \| < \epsilon_1$, all roots of the three polynomials move in a
compact set of $\CC$, since the degrees are constant.

If $(a_t^1,\ldots,a_t^r)$ are the roots common to $g_t,h_t$ and
$k_t$, they are also the roots common to $\tilde{g}_t,
\tilde{h}_t$ and $\tilde{k}_t$. We shall now study the common
roots of two (but not three) of these polynomials.

Suppose that $\alpha$ is a root of $g_0$ and $\beta$ a root of
$h_0$, with $\alpha \neq \beta$. By continuity of the roots of a
family of polynomials, for $A-I$ small enough, the corresponding
roots of $\tilde{g}_0$ and $\tilde{h}_0$ remain distinct.
Applying this remark a finite number of times, we see that any
pair of distinct roots of the polynomials remain distinct after
tranformation by $A$, when $A-I$ is small enough.

Consider now a complex number $b$ which is a root of $g_t$ and
$h_t$, but not $k_t$, so that $g_t(b)=h_t(b)=0$ and $k_t(b)=B
\neq 0$.

Then for $\theta \in \CC$ small enough, replace $g_0$ by $g_0 +
\theta k_0$. We see that $g_0+\theta k_0$ and $h_0$ do not any
more have the common root $b$. By the preceeding remark no new
common root has been created.

Note that the same applies if $b$ is one of the $a_0^s$, by which
we mean that $b$ is a root of the three polynomials, of order
$m+n$ for $g_0$ and $h_0$ and of order $m$ for $k_0$, with $m \ge
1$, $n \geq 1$. Indeed, in this case,
\begin{eqnarray*}
g_0(z) &=& (z-b)^m(z-b)^n\tilde{g}(z), \\
k_0(z) &=& (z-b)^m\tilde{k}(z)
\end{eqnarray*}
and
$$
(g_0+\theta k_0)(z)= (z-b)^m((z-b)^n\tilde{g}(z)+\tilde{k}(z)),
$$
the last factor being non zero at $b$.

Thus $b$ is not any more a common root of $g_0+\theta k_0$ and $h_0$,
except for the multiplicity of the root in all three polynomials.

Repeating this argument a finite number of times, we can replace $g_0,h_0$
and $k_0$ by three new polynomials given by
$$
\left( \begin{array}{c}
\tilde{g}_0 \\
\tilde{h}_0 \\
\tilde{k}_0
\end{array}
\right)
 = A
\left( \begin{array}{c}
{g}_0 \\
{h}_0 \\
{k}_0
\end{array}
\right)
$$
in such a way that the only roots common to two of them are in
fact common to all three.

Defining $\tilde{g}_t, \tilde{h}_t$ and $\tilde{k}_t$ by
(\ref{equ:4.2}) for the same matrix $A$, we see that for $\|t\|$
small enough the roots have the same property.

Thus, the common factor $l_t$ of $g_t,h_t$ and $k_t$ is also the
common factor of, for example, $g_t$ and $h_t$. By Lemma 4.2, it
varies smoothly (resp. complex analytically) with $t$.

\medskip

\noindent {\bf Proof of Proposition \ref{prop:4.1}.} Let $f \in
\Hol_{k,r}^* (\CC P^2)$.  Then identifying $S^2$ with $\CC \cup
\{ \infty \}$ by stereographic projection, $f$ can be represented
(at least on $\CC \subset S^2)$ by a map $p : \CC \to \CC^3
\setminus \{0\}$ with $p(z) =(p_0(z),p_1(z),p_2(z))$, \ $(z \in
\CC)$, a triple  of polynomials with no common zeros and with max
(degree $p_0$, degree $p_1$, degree $p_2)=k$. Now since $f$ has a
finite number of ramification points (in fact no more than $r$)
we can choose the pole of the stereographic projection such that
none of them is at $\infty$\ , then let the ramification points
be $z_1,\ldots,z_t \in \CC$ with ramification indices
$k_1,\ldots,k_t$. Note that $\sum_{i=1}^t k_i = r$.

Next, the first associated curve (see \S 2) $f_{(1)}$  is
represented by $q=p \wedge p' : \CC \to \Lambda^2 \CC^3 \equiv
\CC^3$. This is a triple of polynomials $q=(g,h,k)$ and it is
easily seen that, since $f$ is not ramified at infinity, the
maximum of the degrees of $f,g,k$ is equal to $2k-2$. Further
$(g,h,k)$ have a common zero at $z_i$ of order $k_i$ if and only
if $z_i$ is a ramification point of $f$ of ramification index
$k_i$ so that the highest common factor of $(g,h,k)$ is the {\em
ramification divisor} of $f$ given by
$R(z)=\displaystyle{\prod_{i=1}^t(z-z_i)^{k_i}}$, a polynomial of
degree $r$. Then $q/R$ is a triple of polynomials with no common
roots and, for $z \in \CC$, $f_{(1)}(z)$ is the 2-plane spanned
by $q(z)/R(z)$.

Now suppose that $f_t \in \Hol_{k,r}^*(\CC P^2)$ is a family of
holomorphic maps depending smoothly on a parameter $t \in U
\subset \RR^N$. Then we can choose a family of polynomial maps
$p_t : \CC \to \CC^3 \setminus \{0\}$ as above representing $f_t$
and depending smoothly on $t$ with no ramification point at
infinity.

Here we use the fact that the ramification points, being roots of
polynomials, vary continuously with $t$.

Since the total ramification stays constant ($=r$), the
ramification divisor $R_t$ of $f_t$ has constant degree and we
can apply Lemma 4.5 to see that $R_t$ depends smoothly on $t$.
Hence the corresponding $q_t/R_t$ depends smoothly on $t$ and so
does $(f_t)_{(1)}$. Since $\varphi_t=f_t^\perp \cap f_{t(1)}$ it
is clear that this too varies smoothly with $t$ and the
proposition is proven.

\noindent {\bf Proof of Theorem \ref{th:1.3}.} Similarly, Lemma
\ref{lem:4.5} allows us to conclude that the passage from $f$ to
its conjugate polar is a complex analytic map from
$\Hol_{k,r}^*(\CC P^2)$ to $\Hol_{k',r'}^*(\CC P^2)$, and that
the same applies to its inverse.

\section{The diffeomorphic nature of $G'_{k,r}$}
In this section, we complete the proof of Theorem 1.1. First we
show

\begin{lemma} \label{lem:5.1}  For $k \geq 2$, \ $0 \leq r \leq
\displaystyle{\frac{k+1}{2}}$ and $k \geq 3$,\
$\displaystyle{\frac{4k-11}{3} \leq r \leq \frac{3}{2}k-3}$ \
$G'_{k,r} : \Hol_{k,r}^* (\CC P^2) \to C^j(S^2,{\bf C}P^2)$ has
injective differential at all points $f_0 \in \Hol_{k,r}^* (\CC
P^2)$, for any $j \geq 2$. \end{lemma}

\noindent
{\bf Proof.} Let $f_t \in \Hol_{k,r}^* (\CC P^2)$ be a family of
holomorphic maps depending smoothly on a real parameter $t \in (-
\epsilon, \epsilon)$, $\epsilon > 0$, i.e. a smooth curve in
$\Hol_{k,r}^* (\CC P^2)$. Then, by Proposition \ref{prop:4.1} and
the proof of Theorem \ref{th:1.3}, $\varphi_t = G'_{k,r}(f_t)$
and $g_t$ (the polar of $f_t$) are smooth curves in $C^j(S^2,\CC
P^2)$. Working on a coordinate domain $(U,z)$, let $F_t,
\Phi_t,G_t : U \to \CC^3 \setminus \{0\}$ be families of smooth
maps representing $f_t, \phi_t, g_t$ respectively with $F_t$
holomorphic and $G_t$ antiholomorphic. Then for each $z \in U, \
t \in (-\epsilon, \epsilon)$, $\{F_t(z), \Phi_t(z),G_t(z)\}$ is a
Hermitian orthogonal basis of $\CC^3$. Now $\displaystyle{\frac{d
\varphi_t}{dt} = d \pi \left( \frac{d \Phi_t}{dt} \right)}$
(where, as before, $\pi : \CC^3 \setminus \{0\} \to \CC P^2$ is
the canonical projection). Suppose that $\displaystyle{\frac{d
\varphi_t}{dt}=0}$ for some value of $t$. Then
$\displaystyle{\frac{d \Phi_t}{dt}}$ must be in direction
$\Phi_t$, so that, in particular,
$$ <\frac{d \Phi_t}{dt},F_t>=0.
$$
(Here $<,>$ denotes the standard Hermitian inner product on
${\bf C}^3$). This last equation is equivalent to
$$
<\frac{d F_t}{dt},\Phi_t>=0,
$$
hence
\begin{equation}\label{equ:5.1}
\frac{dF_t}{dt}=\alpha F_t + \beta G_t
\end{equation}
for some smooth functions $\alpha, \beta$ on $U \times
(-\epsilon,\epsilon)$. Differentiating with respect to
$\overline{z}$, since $F_t$ is holomorphic we obtain
$$
0=\frac{\partial \alpha}{\partial \overline{z}}F_t +
\frac{\partial \beta}{\partial \overline{z}}G_t +
\beta\frac{\partial G_t}{\partial \overline{z}}.
$$
Now the triple $\big\{ F_t,\displaystyle{\frac{\partial
G_t}{\partial \overline{z}},G_t} \big\}$ is linearly independent
except at the isolated points where $h_t=\bar{g}_t$ is ramified.
Hence $\beta \equiv 0$ and so from (\ref{equ:5.1})
$$
\frac{d F_t}{dt} = \alpha F_t
$$
which implies that $df_t/dt=0$. The lemma follows.

The proof of Theorem \ref{th:1.1} is completed by

\begin{proposition} \label{prop:5.2} For $k \geq 2$, \ $0 \leq r
\leq \displaystyle{\frac{k+1}{2}}$ and $k \geq 3$,\
$\displaystyle{\frac{4k-11}{3} \leq r \leq \frac{3}{2}k-3}$ \
$G'_{k,r} : \Hol_{k,r}^*(\CC P^2) \to C^j(S^2,\CC P^2)$ is an
embedding with image the closed submanifold
$\Harm_{k-2-r,3k-2-r}(\CC P^2)$ of $C^j(S^2,\CC P^2)$ for any $j
\geq 2$.
\end{proposition}

\noindent
{\bf Proof.} Since $\Hol_{k,r}^*(\CC P^2)$ is finite dimensional,
the differential of $G'_{k,r}$ splits at each point and so by
\cite{Lan} $G'_{k,r}$ is an immersion. By Proposition
\ref{prop:2.5} it is injective and has image
$\Harm_{k-2-r,3k-2-r}(\CC P^2)$, and we show below that it is a
homeomorphism onto its image. Thus it is an embedding and its image
is a closed submanifold of $C^j(S^2,\CC P^2)$. The dimension of
the space $\Harm_{k-2-r,3k-2-r}(\CC P^2)$ is thus equal to the
(real) dimension of $\Hol_{k,r}^* (\CC P^2)$ which is $6k-2r+4$;
the stated dimensions follow easily.

To show that $G'_{k,r}$ is a homeomorphism it suffices to show that
it is proper. To do this, following \cite{Cra},\cite[Lemma
3.3]{F-G-K-O} consider a sequence $(\phi_n)$ which converges to
$\phi$ in $G'_{k,r}(\Hol^*_{k,r}(\CC P^2)) = \Harm_{k-2-r,
3k-2-r}(\CC P^2)$.  Let $f_n = (G'_{k,r})^{-1}(\phi_n) \in
\Hol^*_{k,r}(\CC P^2)$.  It suffices to prove that a subsequence
converges in $\Hol^*_{k,r}(\CC P^2)$.  Note first that
$\Hol^*_{k,r}(\CC P^2)$ can be injected into the projective space
$P(\CC^{3k+3})$ by the map $i:f \mapsto$ projective class of the
coefficients of the polynomials describing $f$.  Since that space
is compact, a subsequence of $(i(f_n))$ converges to $[p] =
[p_0,p_1,p_2]$. (Note that $(p_0, p_1, p_2)$ need not be coprime
and that their maximal degree need not be $k$.)  We retain the
notation $i(f_n)$ for the subsequence.  We can then write $[p] =
[b q_0,b q_1, b q_2]$ where the $q_i$'s are coprime and $[q] =
[q_0, q_1, q_2]$  lies in $\Hol_{k-m}(\CC P^2)$, the space of not
necessarily full holomorphic maps from $S^2$ to $\CC P^2$ of
degree $k-m$, for some $m \geq 0$ . For each $n$, let $h_n$ be the
conjugate polar of $f_n$ belonging to $\Hol_{2k-2-r}^*(\CC P^2)$,
and consider its image by $i$ in the appropriate projective space
of coefficients $P(\CC^{6k-3-3r})$. Again, a subsequence
converges to $[t] = [t_0,t_1,t_2]$, and we have $[t_0,t_1,t_2] =
[a s_0, a s_1, a s_2]$ with the $s_i$'s coprime and of possibly
lower degree than $2k-2-r$.  Now $f_n \perp \overline{g_n}$ for
all $n$ \ so that $q \perp \bar{s}$ and $\psi = (q \oplus
\bar{s})^{\perp}$ is well-defined. Further on $S^2 \setminus
\{\{\mbox{zeros of } a \} \cup \{ \mbox{zeros of } b \} \cup
\{\infty\}\}$, \ $G'(f_n)(x)$ converges to $\psi(x)$ so that
$\phi$ coincides with $\psi$ on a dense set, and therefore
everywhere. Then $3k-2-r = E(\phi) = E(\psi) = \deg q + \deg s$
with $\deg q \leq k$ and $\deg s \leq 2k-2-r$ so that both these
inequalities must be equalities and there can be no loss of
degree above. Thus $q$ must have degree $k$ and ramification
index $r$.  Further, since $r \leq 3k/2 -3$, \ $q$ is full,
otherwise by the Riemann-Hurwitz formula we would have $r = 2k-2
> 3k/2 -3$.  Hence $q \in \Hol_{k,r}^*(\CC P^2)$ as required.

\end{document}